\documentclass{article}
\usepackage{latexsym}

\begin{document}
\title{Explicit formula for even-index Bernoulli numbers}

\author{Dr. Renaat Van Malderen}

\maketitle
\begin{abstract}Bernoulli numbers are usually expressed in terms of their lower index
numbers (recursive). This paper gives an explicit formula for Bernoulli
numbers of even index. The formula contains a remarkable sequence of
determinants.
\end{abstract}

\vspace{10pt}
\section{Explicit expression for even-index Bernoulli numbers.}
\label{section1}

First a sequence of determinants D$_{k}$ (with k= 1, 2, 3, ...) is
defined as follows:

\vspace{10pt}

D$_{1}=\frac{1}{3!}$ \ \ \ \ \ D$_{2 = }\left| {{\begin{array}{*{20}c} {\frac{1}{3!}} \hfill & {\frac{1}{5!}} \hfill \\
 1 \hfill & {\frac{1}{3!}} \hfill \\
\end{array} }} \right|$ \ \ \ \ \ D$_{3}=\left| {{\begin{array}{*{20}c}
 {\frac{1}{3!}} \hfill & {\frac{1}{5!}} \hfill & {\frac{1}{7!}} \hfill \\
 1 \hfill & {\frac{1}{3!}} \hfill & {\frac{1}{5!}} \hfill \\
 0 \hfill & 1 \hfill & {\frac{1}{3!}} \hfill \\
\end{array} }} \right|$

\vspace{10pt}
D$_{4}=\left| {{\begin{array}{*{20}c}
 {\frac{1}{3!}} \hfill & {\frac{1}{5!}} \hfill & {\frac{1}{7!}} \hfill & 
{\frac{1}{9!}} \hfill \\
 1 \hfill & {\frac{1}{3!}} \hfill & {\frac{1}{5!}} \hfill & {\frac{1}{7!}} 
\hfill \\
 0 \hfill & 1 \hfill & {\frac{1}{3!}} \hfill & {\frac{1}{5!}} \hfill \\
 0 \hfill & 0 \hfill & 1 \hfill & {\frac{1}{3!}} \hfill \\
\end{array} }} \right|$ \ \ \ \ \ etc{\ldots}{\ldots}

\begin{flushright}(1.1)\end{flushright}

\begin{flushleft}
Let B$_{2p}$ (with p= 1, 2, 3, {\ldots}) be the even Bernoulli number of 
index 2p.

\vspace{10pt}

B$_{2p}$ is then given by the following formula:

\vspace{10pt}

\[B_{2p} =-2p+\left( {\frac{3}{2}} \right)^{2p}\left\{ 
{1+(2p)!\sum\limits_{k=1}^p {\frac{(-1)^kD_k }{3^{2k}\left[ {2(p-k)} 
\right]!}} } \right\}
\]
\begin{flushright}(1.2)\end{flushright}

in which the D$_{k}$'s are as given above and 0 ! = 1 by definition.
\end{flushleft}

\newpage
\section{Relation between the determinants D$_{k}$.}
\label{section2}
\begin{flushleft}
The determinants D$_{k}$ in (1.1) are in explicit form. Nevertheless a 
recursive relation among the D$_{k}$'s may be given as is obvious from their 
definition (1.1). Putting for simplicity D$_{0}$ = 1, the relation between 
the D$_{k}$'s may be written successively as:

\vspace{10pt}

D$_{1}=\frac{1}{3!}$ D$_{0}$

\vspace{10pt}

D$_{2}=\frac{1}{3!}$ D$_{1 }-\frac{1}{5!}$ D$_{0}$

\vspace{10pt}

D$_{3}=\frac{1}{3!}$ D$_{2 }-\frac{1}{5!}$ D$_{1}+\frac{1}{7!}$ 
D$_{0}$

\vspace{10pt}

D$_{4}=\frac{1}{3!}$ D$_{3 }-\frac{1}{5!}$ D$_{2}+\frac{1}{7!}$ 
D$_{1}-\frac{1}{9!}$ D$_{0}$

\vspace{10pt}

or in general:

\vspace{10pt}

\[
D_k =\sum\limits_{l=1}^k {\frac{(-1)^{l-1}D_{k-l} }{(2l+1)!}} 
\]
\begin{flushright}(2.1)\end{flushright}
\vspace{30pt}

Readers interested in the derivation of (1.2) may contact the author via e-mail.

Dr Renaat Van Malderen

Address: Maxlaan 21, 2640 Mortsel, Belgium

Author can be reached by email: hans.vanmalderen@bertholdtech.com
\end{flushleft}
\end{document}